\documentclass{article}
\usepackage{amssymb}

\title{\bf The Number of Finite Groups Whose Element Orders is Given}
\author{ {\bf Ali Reza Moghaddamfar} and {\bf Wujie
Shi}}

\newenvironment{proof}{\noindent {\sc {Proof}}.}{$\square$
\medskip}

\newtheorem{theorem}{Theorem}

\newtheorem{lm}{Lemma}
\newtheorem{ca}{Case}
\newtheorem{st}{Step}

\begin{document}
\maketitle
\begin{abstract}
\noindent The spectrum $\omega(G)$ of a finite group $G$ is the
set of element orders of $G$. If $\Omega$ is a non-empty subset
of the set of natural numbers, $h(\Omega)$ stands for the number
of isomorphism classes of finite groups $G$ with
$\omega(G)=\Omega$ and put $h(G)=h(\omega(G))$. We say that $G$ is
recognizable (by spectrum $\omega(G)$) if $h(G)=1$. The group $G$
is almost recognizable (resp. nonrecognizable) if $1<h(G)<\infty$
(resp. $h(G)=\infty$). In the present paper, we focus our
attention on the projective general linear groups
$\mbox{PGL}(2,p^n)$, where $p=2^\alpha 3^\beta+1$ is a prime,
$\alpha \geq 0, \beta \geq 0$ and $n\geq 1$, and we show that
these groups cannot be almost recognizable, in other words
$h(\mbox{PGL}(2,p^n))\in \{1, \infty\}$. It is also shown that the
projective general linear groups $\mbox{PGL}(2,7)$ and
$\mbox{PGL}(2,9)$ are nonrecognizable. In this paper a computer
program has also been presented  in order to find out the
primitive prime divisors of $a^n-1$.
\end{abstract}
\renewcommand{\baselinestretch}{1}
\def\thefootnote{ \ }
\footnotetext{{\em $2000$ Mathematics Subject Classification}: 20D05}

\section{Introduction}
Throughout the paper, all the groups under consideration are
finite and simple groups are non-Abelian. For a group $G$, we
denote the set of orders of all elements in $G$ by $\omega(G)$
which has been recently called the {\em spectrum} of $G$.
Obviously $\omega(G)$ is a subset of the set $\mathbb {N}$ of
natural numbers, and it is {\em closed} and {\em partially
ordered} by divisibility, hence, it is uniquely determined by
$\mu(G)$, the subset of its maximal elements.

One of the most interesting concepts in Finite Group Theory which
has recently attracted several researchers is the problem of
characterizing finite groups by element orders. Let $\Omega$ be a
non-empty subset of $\mathbb {N}$. Now, we can put forward the
following questions: {\em Is there any group $G$ with
$\omega(G)=\Omega$} ? {\em If the answer is affirmative then how
many non-isomorphic groups exist with the above set of element
orders} ? Certainly, if there exists such a group, $\Omega$ must
contain 1 and furthermore $\Omega$ must be closed and partially
ordered under the divisibility relation. These conditions are
necessary but not sufficient, for example if $\Omega=\{1, 2, 3,
4, 5, 6, 7, 8, 9 \}$, then there does not exist any group $G$ with
$\omega(G)=\Omega$. In fact, R. Brandl and W. J. Shi in \cite{BS}
have classified all groups whose element orders are consecutive
integers and in that paper they have shown that if
$\omega(G)=\{1, 2, 3, \dots, n \}$, for some group $G$, then
$n\leq 8$.

For a set $\Omega$ of natural numbers, define $h(\Omega)$ to be
the number of isomorphism classes of groups $G$ such that
$\omega(G)=\Omega$, and put $h(G)=h(\omega(G))$. Evidently,
$h(G)\geq 1$. Now we give a {\em ``new classification"} for groups
using the $h$ function. A group $G$ is called {\em recognizable}
(resp. {\em almost recognizable} or {\em nonrecognizable}) if
$h(G)=1$ (resp. $1<h(G)<\infty$ or $h(G)=\infty$). Some list of
simple groups that are presently known to be recognizable, almost
recognizable or nonrecognizable is given in \cite{m4}. In
particular, it was previously known that the projective general
linear groups $\mbox{PGL}(2,2^n)$ with $n\geq 2$ are recognizable
and $\mbox{PGL}(2,2)\cong S_3$ is nonrecognizable (see \cite{S},
Theorem 2).  In \cite{maz}, V. D. Mazurov proved the following
result: Let $P$ be a field which is the union of an ascending
series of finite fields of orders $2^{m_i}, m_i>1, i\in \mathbb{
N}$. If there exists a natural number $s$ such that $2^s$ does
not divide $m_i$ for any $i\in \mathbb{ N}$ then
$h(\mbox{PGL}(2,P))=1$. In all other cases
$h(\mbox{PGL}(2,P))=\infty$. Also he proved the following result
in \cite{m3}: If $p$, $r$ are odd primes, $p-1$ is divisible by
$r$ but not by $r^2$, and $s$ is a natural number non-divisible
by $r$, then $h(\mbox{PGL}(r,p^s))=\infty$.

Let $q=p^n$ where $p$ is a prime. In this paper, we focus our
attention on the projective general linear groups
$\mbox{PGL}(2,q)$. The structure of $\mbox{Aut}(L_2(q))$ is well
known, it is isomorphic to the semidirect product of
$\mbox{PGL}(2,q)$ by a cyclic group of order $n$. On the other
hand we know $\mu (L_2(q))=\{\frac{q-1}{\epsilon}, p,
\frac{q+1}{\epsilon} \}$, $\epsilon =(2,q-1)$, and
$\mu(\mbox{PGL}(2,q))=\{q-1, p, q+1\}$.

A group $G$ is called $C_{pp}$-{\em group} if $p$ is a prime
divisor of $|G|$ and the centralizer of any non-trivial
$p$-element in $G$ is a $p$-group. Evidently, the projective
general linear groups $\mbox{PGL}(2,q)$ where $q=p^n$, are
$C_{pp}$-groups. In \cite{Cpp}, the second author has classified
the simple $C_{pp}$-groups, where $p$ is prime and
$p=2^\alpha3^\beta+1$, $\alpha\geq 0, \beta\geq 0$ (see Lemma
\ref{pp} and Table 1). Using these results, we prove the
following theorem.

\begin{theorem}\label{th1}
Let $p=2^\alpha 3^\beta+1 \ (\alpha \geq 0, \beta \geq 0)$ be a
prime. Then the projective general linear groups
$\mbox{PGL}(2,p^n)$ cannot be almost recognizable. In other
words, $h(\mbox{PGL}(2,p^n))\in \{1, \infty\}$.
\end{theorem}

In 1994, R. Brandl and W. J. Shi in \cite{BSH} showed that all
projective special linear groups $L_2(q)$ with $q\neq 9$ are
recognizable and $L_2(9)$ is nonrecognizable.

Here, we similarly prove that:

\begin{theorem}\label{th2}
The projective general linear group $\mbox{PGL}(2,9)$ is
nonrecognizable.
\end{theorem}

For us it was interesting to face with some groups $G$ such that
$\mu(G)$ contain three consecutive natural numbers in the form
$\{p-1, p, p+1\}$ where $p\geq 5$ is a prime. Such sets appear
for almost simple groups $\mbox{PGL}(2,p)$, where $p\geq 5$ is a
prime, in fact we proved in \cite{MS} that $h(\mbox{PGL}(2,p))\in
\{1, \infty \}$. For $\infty$ we have found an example. It has
been proved in \cite{BS} that $\mbox{PGL}(2,5)\cong S_5$ has
$\infty$ for its $h$ function. Here we also give another example
of groups of type $\mbox{PGL}(2,p)$ with value $\infty$ for its
$h$ function.
\begin{theorem}\label{th3}
There exists an extension $G$ of a $7$-group by $2.S_4$ such that
$\mu(G)=\mu(\mbox{PGL}(2,7))=\{6,7,8\}$. In particular, the
projective general linear group $\mbox{PGL}(2,7)$ is
nonrecognizable.
\end{theorem}
{\em Notation.} Our notation and terminology are standard (see
\cite{atlas}). Given a group $G$, denote by $\pi(G)$ the set of
all prime divisors of the order of $G$. If $m$ and $n$ are natural
numbers and $p$ is a prime, then we let $\pi(n)$ be the set of
all primes dividing $n$, and $r_{[n]}$ the largest prime not
exceeding $n$. Note that $\pi(G)=\pi(|G|)$. The notation
$p^m\parallel n$ means that $p^m|n$ and $p^{m+1}\nmid n$. The
expression $G=K:C$ denotes the split extension of a normal
subgroup $K$ of $G$ by a complement $C$.

\section{Some Preliminary Results}
First, we collect some results from Elementary Number Theory which
will be useful tools for our further investigations in this
paper. We start with a famous theorem due to Zsigmondy. \\[0.4cm]
{\bf Zsigmondy's Theorem} (see \cite{zsigmondy}). {\em Let $a$
and $n$ be integers greater than $1$. Then there exists a
``primitive prime divisor" of $a^n -1$, that is a prime $s$
dividing $a^n -1$ and not dividing $a^i -1$ for $1\leq i\leq n-1$,
except if

$(1)$ $a=2$ and $n=6$, or

$(2)$ $a$ is a Merssene prime
and $n=2$. }\\[0.4cm]
\indent  We denote by $a_n$ one of these primitive prime divisors
of $a^n -1$. Evidently, if $a_n$ is a primitive prime divisor of
$a^n-1$, then $a$ has order $n$ modulo $a_n$ and so $a_n \equiv
1\!\!\!\!\pmod{n}$. Thus $a_n\geq n+1$.

The next elementary result will be needed later.
\begin{lm}\label{prim}
Let $p$ and $q$ be two primes and $m$ be a natural number, where
$p, q$ and $m$ satisfying one of the following conditions. Then,
for every $n\geq m$, there exists a primitive prime divisor
$p_n>q$.

\begin{tabular}{lllll}
$(1)$ & $p=7$,& $m=5$& and &$q=13$,\\

$(2)$ & $p=13$,& $m=5$& and &$q=19$,\\

$(3)$ & $p=17$,& $m=4$& and &$q=19$,\\

$(4)$ & $p=19$,& $m=7$& and &$q=37$,\\

$(5)$ & $p=37$, & $m=7$& and &$q=109$,\\

$(6)$ & $p=73$,& $m=5$& and &$q=127$.\\
\end{tabular}
\end{lm}
\begin{proof} In all cases, if $n\leq q$, the result is
straightforward. Therefore, we may assume that $n>q$. Since
$$\pi(q!)\subseteq \pi(p\prod_{i=1}^{q-1}(p^i-1))\subset
\pi(p\prod_{i=1}^{n}(p^i-1)),$$ by Zsigmondy's theorem we deduce
that there exists a primitive prime divisor $p_n>q$, completing
the proof. \end{proof}

\noindent {\sc Function for finding the primitive prime divisors}.
\ In the following we submit a \textsf{GAP} program \cite{gap},
which
determines all the primitive prime divisors in the sequence $a^i-1$ ($i=1,2, \ldots, n$) for some $a$ and $n$.\\[0.5cm]
{\tt gap>  PPD:=function(a,n)\\
\hspace*{0.6cm} local  b,i,j,s1,s2,s;\\
\hspace*{0.9cm} for i in [1..n] do\\
\hspace*{0.9cm} s1:=Set(Factors(a$^{\wedge}$i-1));\\
\hspace*{0.6cm} s2:=[];\\
\hspace*{1.5 cm} for j in [1..(i-1)] do\\
\hspace*{2.2cm} b:=Set(Factors(a$^{\wedge}$j-1));\\
\hspace*{2.2cm} Append(s2,b);\\
\hspace*{1.5cm}       od;\\
\hspace*{1.5cm}            s:=Difference(s1,s2);\\
\hspace*{1.5cm}            Print(i,"     ",s,"$\backslash$n");\\
\hspace*{0.6cm}      od;\\ end;\\ }

Using this programme we list all primitive prime divisors $p_n$
for $p=$ 7, 13, 17 and $2\leq n\leq 19$, in Table 1. Using Table
1, the reader can easily check the proof of Lemma \ref{prim}
(1)-(3) for $n\leq q$.

\begin{lm}\label{ex}
Let $p$ and $q$ be two primes and $m, n$ be natural numbers such
that $p^m=q^n+1.$ Then one of the following holds:

$(1)$ $n=1$, $m$ is a prime number, $p=2$ and $q=2^m-1$ is a
Mersenne prime;

$(2)$ $m=1$, $n$ is a power of $2$, $q=2$ and $p=2^n+1$ is a
Fermat prime;

$(3)$ $p=n=3$ and $q=m=2$.
\end{lm}
\begin{proof}
Well known exercise using the Zsigmondy's theorem.
\end{proof}

The set $\omega(G)$ defines the {\em prime graph} $\mbox{GK} (G)$
of $G$ whose vertex-set is $\pi(G)$ and two primes $p$ and $q$ in
$\pi(G)$ are adjacent (we write $p\sim q$) if and only if $pq\in
\omega(G)$. The number of connected components of $\mbox{GK}(G)$
is denoted by $t(G)$, and the connected components are denoted by
$\pi_i=\pi_i(G)$, $i=1,2,\ldots,t(G)$. If $2\in \pi(G)$ we always
assume $2\in \pi_1$. Denote by $\mu_i(G)$ the set of all $n\in
\mu(G)$ such that $\pi(n)\subseteq \pi_i$.\\[0.5cm]
\noindent {\bf The Gruenberg-Kegel Theorem} (see \cite{williams}).
{\em If $G$ is a group with disconnected graph $\mbox{GK}(G)$
then one of the following holds:

$(1)$ $t(G)=2$, $G$ is Frobenius or $2$-Frobenius.

$(2)$ $G$ is an extension of a $\pi_1(G)$-group $N$ by a group
$G_1$, where $S\leq G_1 \leq \mbox{Aut}(S)$, $S$ is a simple
group and $G_1/S$ is a $\pi_1(G)$-group. Moreover $t(S)\geq t(G)$
and for every $i$, $2\leq i\leq t(G)$, there exists $j$, $2\leq
j\leq t(S)$ such that $\mu_j(S)=\mu_i(G)$. }

\begin{center}
\small {\bf Table 1} {\em The primitive prime divisors $p_n$ where
$p\in \{7,13,17\}$ and $2\leq n\leq 19$.}\\[0.5cm]

\begin{tabular}{|l|l|l|l|}
\hline $n$ & $7_n $ & $13_n $ & $17_n $
\\
\hline 2&--&7&3\\
\hline 3&19&61&307\\
\hline 4&5&5, 17&5, 29\\
\hline 5&2801&30941&88741\\
\hline 6&43& 157&7, 13\\
\hline 7&29, 4733&5229043&25646167\\
\hline 8&1201&14281&41761\\
\hline 9 &37, 1063&1609669&19, 1270657\\
\hline  10&11, 191&11, 2411&11, 71, 101\\
\hline 11&1123,& 23, 419, &2141993519227\\
         &293459&859, 18041&\\
\hline 12&13, 181&28393&83233\\
\hline 13&16148168401&53, 264031&212057,\\
         &&1803647&2919196853\\
\hline 14&113, 911&29, 22079&22796593\\
\hline 15&31, 159871&4651, 161971&6566760001\\
\hline 16&17, 169553&407865361&18913,\\
          &&&184417\\
\hline 17&14009&103, 443,&10949,\\
         &2767631689&15798461357509&1749233,\\
         &&&2699538733\\
\hline 18&117307&19, 271,&1423, \\
         &&937&5653\\
\hline 19&419&12865927,&229, 1103, \\
          &4534166740403&9468940004449&202607147,\\
& & &291973723\\

\hline
\end{tabular}
\end{center}
\vspace{0.2 cm}

\begin{lm} \label{mu} Let $S$ be a simple group
with disconnected prime graph $\mbox{GK}(S)$. Then $|\mu_i(S)|=1$
for $2\leq i\leq t(S)$. Let $n_i(S)$ be a unique element of
$\mu_i(S)$ for $i\geq 2$. Then value for $S$, $\pi_1(S)$ and
$n_i(S)$ for $2\leq i\leq t(S)$ are the same as in Tables $2a-2c$
of \cite{m4}.
\end{lm}
\begin{proof}
The simple groups $S$ and the sets of $\pi_i(S)$ are described in
\cite{williams} and \cite{kondratev}; the rest is proved in Lemma
4 of \cite{km}. The values of the numbers $n_i(S)$, $i\geq 2$ are
listed in Table 2a-2c of \cite{m4}.
\end{proof}

We also use the following lemma (see \cite{m3}, Lemma 1).

\begin{lm}\label{Maz1} If a group $G$ contains a
soluble minimal normal subgroup then $G$ is nonrecognizable. In
particular, if $G$ is a soluble group then $G$ is nonrecognizable.
\end{lm}

The following result of V. D. Mazurov will be used several times.

\begin{lm}(see \cite{m2}) \label{Maz3}
Let $G$ be a group, $N$ a normal subgroup of $G$, and $G/N$ a
Frobenius group with Frobenius kernel $F$ and cyclic complement
$C$. If $(|F|,|N|)=1$ and $F$ is not contained in $NC_{G}(N)/N$,
then $p|C|\in \omega(G)$ for some prime divisor $p$ of $|N|$.
\end{lm}

The following lemma is taken from (\cite{S}, Theorem 2).

\begin{lm}\label{item} Let $G$ be a group such that
$$\mu(G)=\mu(\mbox{PGL}(2,2^n))=\{2^n-1,2, 2^n+1\}.$$ Then, the following statements hold.

$(1)$ If $n\geq 2$, then $G \cong \mbox{PGL}(2,2^n)$.

$(2)$ If $n=1$, then $G \cong S_3$ has $\infty$ for its $h$
function.
\end{lm}

We are now ready to prove the following lemma.

\begin{lm}\label{item} Let $G$ be a group such that
$$\mu(G)=\mu(\mbox{PGL}(2,p^n))=\{p^n-1,p, p^n+1\},$$ where $p$ is an odd
prime, $n\geq 2$. Then, the following statements hold.

$(1)$ If $(p,n)\neq (3,2)$, then item $(2)$ of the
Gruenberg-Kegel theorem holds. Moreover, $S$ is isomorphic to none
of the following simple groups:

$(a)$ alternating groups on $n\geq 5$ letters,

$(b)$ sporadic simple groups,

$(c)$ $L_2(p^k)$ where $k\neq n$, or

$(d)$ $L_2(2p^m\pm 1)$, $m\geq 1$, where $2p^m\pm 1$ is a prime.

$(2)$ If $(p,n)=(3,2)$, then there exists a soluble group $G$
such that $\mu(G)=\mu(\mbox{PGL}(2,3^2))$.
\end{lm}

\begin{proof} (1) First of all, we show that $G$ is insoluble.
Assume the contrary. If $\pi(p^n-1)=\{2\}$, then by Lemma
\ref{ex} we obtain $(p,n)=(3,2)$ which is a contradiction. Hence,
there exists a prime $2\neq r\in \pi(p^n-1)$. On the other hand,
we consider the primitive prime divisor $s=p_{2n}$. Now assume
that $H$ is a $\{p,r,s\}$-Hall subgroup of $G$. Since $G$ has no
elements of order $pr, ps$ and $rs$, it follows that $H$ is a
soluble group all of whose elements are of prime power orders. By
(\cite{Higman}, Theorem 1), we must have $|\pi(H)| \leq 2$, which
is a contradiction.

Since $t(G)=2$, $G$ satisfies the conditions of the
Gruenberg-Kegel theorem. Now we show that $G$ is neither Frobenius
nor 2-Frobenius. Evidently, $G$ can not be a 2-Frobenius group,
because $G$ is insoluble. Suppose $G=KC$ is a Frobenius group
with kernel $K$ and complement $C$. Clearly $C$ is insoluble,
$\pi(C)=\pi_1(G)=\pi(p^{2n}-1)$, $\pi(K)=\pi_2(G)=\{p\}$ and by
(\cite{pass}, Theorem 18.6) $C$ has a normal subgroup $C_0$ of
index $\leq 2$ such that $C_0\cong SL(2,5)\times Z$, where every
Sylow subgroup of $Z$ is cyclic and $\pi(Z)\cap
\pi(30)=\emptyset$. Therefore $\mbox{GK}(C)$ can be obtained from
the complete graph on $\pi(C)$ by deleting the edge $\{3,5\}$. On
the other hand, if there exist primes $2\neq r\in \pi(p^n-1)$ and
$2\neq s\in \pi(p^n+1)$, then since $rs\notin \omega(G)$ it
follows that $rs\notin \omega(C)$. Hence, we must have $Z=1$ and
$\pi(p^{2n}-1)=\pi(SL_2(5))=\{2, 3, 5\}$ and since $\{2, 3,
5\}\subset \pi(p^4-1)$, by Zsigmondy's theorem we obtain that
$n=2$. Now, it is easy to see that $\pi(p^2-1)=\{2,3\}$ and
$\pi(p^2+1)=\{2,5\}$. From $\pi(p^2-1)=\{2,3\}$, we infer that
$p$ is a Mersenne prime or a Fermat prime. In the first case we
obtain $p=7$, and in the latter case $p=17$. If $p=17$, then $29
\in \pi(p^2+1)$, a contradiction. If $p=7$, then $C$ contains an
element of order $16$, which is a contradiction.

Therefore, by the Gruenberg-Kegel theorem, $G$ is an extension of
a $\pi_1(G)$-group $N$ by a group $G_1$, where $S\leq G_1 \leq
\mbox{Aut}(S)$, $S$ is a simple group and $G_1/S$ is a
$\pi_1(G)$-group. Now, we show that $S$ is not isomorphic to an
alternating group, a sporadic simple group, a linear group
$L_2(p^k)$ where $k\neq n$ or $L_2(2p^m\pm 1)$, $m\geq 1$, where
$2p^m\pm 1$ is a prime.

Before beginning we recall that in the prime graph of $G$ the
connected component $\pi_1(G)$ consists of the primes in
$\pi(p^n-1)$ which form a complete subsection and also the primes
in $\pi(p^n+1)$ which forms another complete subsection.
Moreover, every odd vertex in $\pi(p^n-1)$ is not joined to any
odd vertex in $\pi(p^n+1)$.

(a) \ Assume that $S \cong A_m$, $m \geq 5$. By Lemma \ref{mu},
$m=p, p+1, p+2$. Suppose $S\cong A_p$, $p \geq 5$. We have that
in the prime graph $\mbox{GK}(A_p)$ the vertex 3 is joined to $2,
5, 7, \ldots, r_{[p-3]}$.  If 3 divides $p-1$, then by the remark
mentioned in the previous paragraph, we conclude that $2, 3, 5,
\ldots, r_{[p-3]}$ belong to $\pi(p^n-1)$. Now, if there exists a
prime $s \in \pi (p^n+1)\backslash \pi (A_p)$ then $s \in \pi(N)$,
because $ A_p \cong S \leq G/N \leq \mbox{Aut}(S)\cong S_p$. On
the other hand, $A_4=2^2:3 \leq A_p$ and by lemma \ref{Maz3} it
follows that $ s \sim 3$ which is a contradiction. Hence, $
\pi(p^n+1) \subseteq \pi(A_p)$. As $(p^n-1, p^n+1)=2$ and $2, 3,
5, 7, \ldots, r_{[p-3]} \in \pi(p^n-1)$, the only possible cases
are: $\pi(p^n+1)=\{2\}$ or $\pi(p^n+1)=\{2,p-2\}$ in which in the
latter case $p-2$ is a prime. Evidently, the first case will
never occur. So, we consider the case $\pi(p^n+1)=\{2,p-2\}$,
i.e., $p^n+1=2^l(p-2)^k$. Now, if $k>1$ then since $(p-2)^k
\in\omega(G)$ and $(p-2)^k \notin
\omega(\mbox{Aut}(S))=\omega(S_p)$ we obtain $(p-2)\in\pi(N)$ and
again since $A_p$ contains a Frobenius subgroup of shape $2^2:3$
by Lemma \ref{Maz3}, we get $p-2\sim 3$ which is a contradiction.
Finally, we have $k=1$ and $l>1$. Moreover $2\parallel p^n-1$
which implies that $n$ must be odd. But in this case we have
$p^n+1=(p+1)(p^{n-1}-p^{n-2}+\cdots-p+1)=2^l(p-2)$ for which it
follows that $p^{n-1}-p^{n-2}+\cdots-p+1=p-2$, giving no solution
for $p\geq5$. This final contradiction shows that $S\not\cong
A_p$. The case when 3 divides $p+1$, is similar. The other cases
are settled similarly.

(b) \ Suppose $S$ is isomorphic to one of the sporadic simple
groups, for instance $S\cong J_2$. Since $p\in \pi_2(G)$, by
Lemma \ref{mu} it follows that $p=7$. If $n\geq 5$, then we
choose the primitive prime divisors $7_n, 7_{2n}$ in $\pi(G)$.
Evidently, $7_{2n}\in \pi(p^n+1)$, and so $G$ does not contain an
element of order $7_n.7_{2n}$. On the other hand since
$\pi(\mbox{Aut}(S))=\{2,3,5,7\}\subset
\pi(7\prod_{i=1}^4(7^i-1))$, it follows that $7_n,7_{2n}\notin
\pi(\mbox{Aut}(S))$. Therefore $7_n, 7_{2n}\in \pi(N)$, and since
$N$ is nilpotent we obtain that $7_n.7_{2n}\in \omega(N)$, which
is a contradiction. Thus $n\leq 4$. If $n=4$, then
$\mu(G)=\{2^5.3.5^2,7,2.1201\}$. Because, there does not exist
any element of order 1201 in $\mbox{Aut}(S)$, 1201 divides the
order of $N$. Without loss of generality we may assume that
$N\neq 1$ is an elementary Abelian 1201-group. Now since $S$
contains the Frobenius group  $A_4=2^2:3$, from Lemma \ref{Maz3}
we infer that $G$ contains an element of order $1201.3$, which is
a contradiction. If $n=2$ or $3$, then $5\in \pi(S)\backslash
\pi(G)$, which is impossible.

The other sporadic simple groups are examined similarly.

(c) \ Assume that $S\cong L_2(p^k)$, where $k\neq n$. In this
case we must have $k<n$, since otherwise by Zsigmondy's Theorem
we get $p_{2k}\in \pi(S)\backslash \pi(G)$, which is a
contradiction. Now, we choose the primitive prime divisors $p_n$
and $p_{2n}$ in $\pi(G)$. Since $p_{2n}\in \pi(p^n+1)$, $G$ does
not contain an element of order $p_n.p_{2n}$. On the other hand,
since $p_n>n>k$ we have $p_n, p_{2n}\notin
\pi(\mbox{Aut}(S))=\pi(\mbox{PGL}(2,p^k)\rtimes Z_k)$, and so
$p_n, p_{2n}\in \pi(N)$. Now, since $N$ is nilpotent we obtain
that $p_n.p_{2n}\in \omega(N)\subset \omega(G)$, which is a
contradiction.

(d) \ Suppose that $S\cong L_2(2p^m\pm 1)$, $m\geq 1$, where
$2p^m\pm 1$ is a prime. By the structure of $\mu(G)$, we see that
$p\in \omega(G)$ and $p^2\not \in \omega(G)$. So, if $S\cong
L_2(2p^m\pm 1)$, where $2p^m\pm 1$ is a prime and $m\geq 1$, then
we deduce $m=1$, because in this case $p^m\in \omega(L_2(2p^m\pm
1))=\omega(S)\subseteq \omega(G)$. On the other hand, we know
$|\mbox{Aut}(S)|=2^2p(p\pm 1)(2p\pm 1)$, where $2p\pm 1$ is a
prime, and so $\pi(\mbox{Aut}(S))=\{p, 2p\pm 1\}\cup \pi(p\pm
1)$. If $n=2$, then $(2p\pm 1,|G|)=1$, which is a contradiction.
Therefore $n\geq 3$. Now, we consider the primitive prime
divisors $p_n$ and $p_{2n}$. Since $(p_n,p_{2n})=(p_n, p\pm
1)=(p_{2n}, p\pm 1)=1$, it follows that $p_n\notin
\pi(\mbox{Aut}(S))$ or $p_{2n}\notin \pi(\mbox{Aut}(S))$, thus we
may assume $N$ is a $p_n$-subgroup or a $p_{2n}$-subgroup. First,
we assume that $2p+1$ is a prime. Let $P$ be a Sylow
$(2p+1)$-subgroup of $S$, then $N_S(P)$, the normalizer of $P$ in
$S$, is a Frobenius group of order $(2p+1)p$, with cyclic
complement of order $p$. Now, by Lemma \ref{Maz3}, we deduce that
$p_n\sim p$ or $p_{2n}\sim p$, which is a contradiction. Next, we
assume that $2p-1$ is a prime. In this case, if there exists a
prime $s \in \pi (p^n+1)\backslash \pi (\mbox{Aut}(S))$ then $s
\in \pi(N)$, because $G/N \leq \mbox{Aut}(S)$. Moreover, if $Q$
is a Sylow $(2p-1)$-subgroup of $S$, then $N_S(Q)$ is a Frobenius
group of order $(2p-1)(p-1)$, with cyclic complement of order
$p-1$. Now, as previous case we get $s.(p-1)\in \omega(G)$, which
is a contradiction. Hence, $ \pi(p^n+1) \subseteq
\pi(\mbox{Aut}(S))$. As $(p^n-1, p^n+1)=2$, the only possible case
is $\pi(p^n+1)=\{2,2p-1\}$, i.e., $p^n+1=2^l(2p-1)^k$ for some $l$
and $k$ in $\mathbb{N}$. Now, if $k>1$ then since $(2p-1)^k
\in\omega(G)$ and $(2p-1)^k \notin \omega(\mbox{Aut}(S))$ we
obtain $(2p-1)\in\pi(N)$. On the other hand, it is easy to see
that $p_n\in \pi(G)\backslash \pi(\mbox{Aut}(S))$, and so $p_n\in
\pi(N)$. Since $N$ is nilpotent, we deduce that $p_n\sim (2p-1)$,
which is a contradiction. Finally, we have $k=1$ and since
$(p,n)\neq (3,2)$, we obtain that $l>1$. Moreover $2\parallel
p^n-1$ which implies that $n$ must be odd. But in this case we
have $p^n+1=(p+1)(p^{n-1}-p^{n-2}+\cdots-p+1)=2^l(2p-1)$ for
which it follows that $p^{n-1}-p^{n-2}+\cdots-p+1=2p-1$, giving
no solution for $p\geq 3$. This final contradiction shows that
$S\not\cong L_2(2p^m\pm 1)$.

  (2) Consider the group $H=\langle a,b|a^8=b^5=1,
ba=ab^2\rangle \cong Z_5:Z_8$. For this group we have
$\mu(H)=\{8,10\}$. Now, we assume that $G$ is an extension of
elementary Abelian 3-group $K$ of order $3^{40l}$ by $H$, and the
generators $a$, $b$ of $H$ act on $K$ cyclically. Then $G$ is a
soluble group and
$\omega(G)=\omega(\mbox{PGL}(2,3^2))=\{1,2,3,4,5,8,10\}.$
\end{proof}

The following lemma gives a classification of simple $C_{pp}$-
groups, where $p$ is a prime of form $p=2^\alpha3^\beta+1$,
$\alpha\geq 0, \beta\geq 0$.
\begin{lm} \label{pp}(see \cite{Cpp}) Let $p$ be a prime and
$p=2^\alpha 3^\beta +1$, $\alpha \geq 0, \beta \geq 0$. Then any
simple $C_{pp}$-group is given by Table $2$.
\end{lm}

The next lemma gives the maximal odd factors set $\psi(F_4(q))$ of
$\mu(F_4(q))$, $q=2^e$.
\begin{lm} \label{F_4(q)} Let $S\cong F_4(q)$,
where $q=2^e$, $e\geq 1$. Then $\psi(S)= \{q^4-1, q^4+1,
q^4-q^2+1,(q-1)(q^3+1), (q+1)(q^3-1)\}$.
\end{lm}
\begin{proof} The $2'$-elements of $S$ is contained in the maximal tori of $S$. From \cite{Shinoda}
we see that $\mu(F_4(q))$ contains 25 maximal tori $H(1)$,
$H(2)$, $\dots$, $H(25)$. Since $(q-1,q^3+1)=1, (q+1,q^3-1)=1$,
$H(13)$ and $H(15)$ are all cyclic. The conclusion holds.
\end{proof}

\section{Main Results}

In this section we prove the statement of Theorems \ref{th1}
, \ref{th2} and \ref{th3}.\\[0.3cm]
\textsc{Proof of Theorem \ref{th1}}.  Let $G$ be a group and
$$\mu(G)=\mu(\mbox{PGL}(2,p^n))=\{p^n-1,p,p^n+1\},$$
where $p=2^\alpha3^\beta+1$ is a prime, and $n$ is a natural
numbers. If $\alpha=\beta=0$, then $p=2$ and the result is correct
by Lemma 6. Also for $n=1$, the result holds by \cite{MS}, and so
from now on we assume that $p$ is an odd prime and $n\geq 2$. Then
$t(G)=2$, in fact we have
\begin{center}
$\pi_1(G)=\pi(p^{2n}-1)$ \ \ \ and \ \ \  $\pi_2(G)=\{p\}.$
\end{center}

Lemma \ref{item}(1) shows that $G$ is an extension of a
$\pi_1(G)$-group $N$ by a group $G_1$, where $S\leq G_1 \leq
\mbox{Aut}(S)$, $S$ is a simple group of Lie type (except
$L_2(p^k)$, $k\neq n$ and $L_2(2p^m\pm 1)$ where $m\geq 1$ and
$2p^m\pm 1$ is a prime) and $G_1/S$ is a $\pi_1(G)$-group.
Moreover, there exists $2\leq j\leq t(S)$ such that
$\mu_j(S)=\{p\}$, in fact $S$ is a simple $C_{pp}$-group. Using
the results summarized in Table 2, we will show that $S$ is
isomorphic to $L_2(p^n)$.

\begin{st}
$S\cong L_2(q)$, $q=p^n$, $n\geq 2$.
\end{st}
In the following case by case analysis we assume that $S\ncong
L_2(p^n)$ and try to obtain a contradiction. Moreover, as $S$ is
always a $C_{pp}$-group for some appropriate prime $p$, we make
use of the results summarized in Table 2 and Lemma 7 and omit the
details of the argument.

\begin{ca}
$q=3^n$, $n\geq 2$.
\end{ca}
In this case $S$ can only be isomorphic to one of the following
simple groups: $L_2(2^3)$, $L_3(2^2)$. Since $G$ does not contain
an element of order $9$, $S$ can not be isomorphic to $L_2(2^3)$.
If $S\cong L_3(2^2)$, then since $7\in \pi(S)$ we obtain that
$n\geq 6$. Assume first that $n=6$. In this case we have
$\pi(G)=\{2, 3, 5, 7, 13, 73\}$. Evidently $13,
73\notin\pi(\mbox{Aut}(S))$ and $13\nsim 73$. Hence $\{13,
73\}\subseteq \pi(N)$, and since $N$ is nilpotent we get
$13.73\in \omega(N)$, which is a contradiction.
 Next we suppose that $n\geq7$.
\begin{center}
{\bf Table 2} {\em Simple $C_{pp}$-groups, $p=2^\alpha
3^\beta+1$, $\alpha \geq 0, \beta \geq 0$.}\\[0.3cm]

\begin{tabular}{|l|l|}
\hline
$p$ &  simple $C_{pp}$-groups\\[0.1cm]
\hline $2$ & $A_5$, $A_6$,
  $L_2(q)$ where $q$ is a Fermat prime, a
Mersenne prime or \\ & $q=2^m$, $m\geq 3$, $L_3(2^2)$,
$Sz(2^{2m+1})$,
$m\geq 1$.\\[0.1cm]
 \hline
$3$ & $A_5$, $A_6$, $L_2(q)$, $q=2^3$, $3^m$ or $2.3^m \pm 1,$
which is a prime, $m\geq 1$,\\ &
  $L_3(2^2)$\\[0.1cm]
\hline $5$ & $A_5$, $A_6$, $A_7$,
 $M_{11}$, $M_{22}$,
 $L_2(q)$, $q=7^2$, $5^m$ or $2.5^m \pm 1$, which is
\\ & a prime,
$m\geq 1,$ $L_3(2^2)$, $S_4(q)$, $q=3, 7$, $U_4(3)$, $Sz(q)$,
$q=2^3, 2^5.$\\[0.1cm]
\hline $7$ & $A_7$, $A_8$, $A_9$, $M_{22}$, $J_1$, $J_2$, $HS$,
$L_2(q)$, $q=2^3$, $7^m$ or $2.7^m-1$,\\ &
 which is a  prime, $m\geq 1$, $L_3(2^2)$, $S_6(2)$, $O_8^+(2)$, $G_2(q)$,
  $q=3, 19$,\\ & $U_3(q)$, $q=3, 5, 19, $
 $U_4(3)$, $U_6(2)$, $Sz(2^3)$.\\[0.1cm]
\hline $13$ & $ A_{13}$, $A_{14}$, $A_{15}$, $Suz$, $Fi_{22}$,
$L_2(q)$, $q=3^3$, $5^2$, $13^m$ or $2.13^m-1$,\\ & which is a
prime,
 $m \geq 1$, $L_3(3)$, $L_4(3)$, $O_7(3)$, $S_4(5)$, $S_6(3)$,
  \\ &  $O_8^+(3)$, $G_2(q)$, $q=2^2, 3$,  $ F_4(2)$,
  $U_3(q)$, $q=2^2, 23$, $Sz(2^3)$,  \\ &
  $^3D_4(2)$, $^2E_6(2)$, $^2F_4(2)'.$\\[0.1cm]
\hline $17$ & $A_{17}$, $A_{18}$, $A_{19}$, $J_3$, $He$,
$Fi_{23}$, $Fi'_{24}$, $L_2(q)$, $q=2^4, 17^m$ or \\ & $2.17^m\pm
1,$ which is a prime $, m\geq 1$, $S_4(4)$, $S_8(2)$, $F_4(2)$,\\
& $O_8^-(2)$, $O_{10}^-(2)$,
$^2E_6(2). $\\[0.1cm]
\hline $19$ & $A_{19}$, $A_{20}$, $A_{21}$, $J_1$, $J_3$, $O'N$,
$Th$, $HN$, $L_2(q)$, $q=19^m$ \\ & or $2.19^m-1,$
  which is a prime, $m \geq 1$,  $L_3(7)$, $U_3(2^3)$,
  \\ & $R(3^3)$, $^2E_6(2).$\\[0.1cm]
\hline $37$ & $A_{37}$, $A_{38}$, $A_{39}$, $J_4$, $Ly$, $L_2(q)$,
$q=37^m$ or $ 2.37^m-1,$ \\ & which is a prime,
$m\geq 1,$ $U_3(11)$, $R(3^3)$, $^2F_4(2^3).$\\[0.1cm]
\hline $73$ & $A_{73}$, $A_{74}$, $A_{75}$, $L_2(q)$, $q=73^m$ or
$2.73^m-1$, which is a prime,\\ & $m\geq 1$, $L_3(2^3)$,
$S_6(2^3)$, $G_2(q)$, $q=2^3$, $3^2$, $F_4(3)$,
$E_6(2)$, $E_7(2)$, \\ & $U_3(3^2)$, $^3D_4(3).$\\[0.1cm]
\hline $109$ & $A_{109}$, $A_{110}$, $A_{111}$, $L_2(q)$,
$q=109^m$ or $2.109^m-1$, which is  \\ &  a prime, $m\geq 1$, $^2F_4(2^3). $\\[0.1cm]
\hline $p=$ & $A_p$, $A_{p+1}$, $A_{p+2}$, $L_2(q)$, $q=2^m, p^k,$
$2p^k\pm 1$, which is a prime,\\
$2^m+1$,& $k\geq 1$, $S_a(2^b)$, $a=2^{c+1}$ and $b=2^d$, $c\geq
1$,
$c+d=s$, $F_4(2^e)$, \\
$m=2^s$ & $e\geq 1$, $4e=2^s$, $ O_{2(m+1)}^-(2)$, $s\geq 2$,
$O_a^-(2^b)$, $a=2^{c+1}$ and $b=2^d$, \\
 & $c\geq 2$, $c+d=s.$\\
\hline Other & $ A_{p}$, $A_{p+1}$, $A_{p+2}$, $L_2(q)$, $q=p^m$
or $2p^m-1,$
 which is a prime,\\ & $m\geq 1. $\\[0.1cm]
\hline
\end{tabular}
\end{center}
\vspace{0.2cm}

Now we choose the primitive prime divisors $3_n$ and $3_{2n}$ in
$\pi(G)$. Evidently $3_{2n}\in \pi(3^n+1)$, and so $3_n\nsim
3_{2n}$. Moreover, since $\{2, 3, 5, 7, 11,
13\}=\pi(3\prod_{i=1}^6(3^i-1))$, $3_n, 3_{2n}\notin
\pi(\mbox{Aut}(S))$, and hence $3_n, 3_{2n}\in \pi(N)$. Again
since $N$ is nilpotent, $N$ contains an element of order $3_n.
3_{2n}$, which is of course impossible.
\begin{ca}
$q=5^n$, $n\geq 2$.
\end{ca}
In this case we see that $S$ can only be isomorphic to one of the
following simple groups: $L_2(7^2)$, $L_3(2^2)$, $S_4(3)$,
$S_4(7)$, $U_4(3)$, $Sz(2^3)$ or $Sz(2^5)$. Since $G$ has no
element of order 25, $S$ can not be isomorphic to $L_2(7^2)$ or
$Sz(2^5)$. If $S$ is isomorphic to one of the simple groups:
$L_3(2^2)$, $S_4(7)$, $U_4(3)$, or $Sz(2^3)$, then $7\in \pi(S)$
and so we must have $n\geq 6$. Also note that
$$ \pi(S)\subset \{2,3,5,7, 11, 13\}\subset \pi(5\prod_{i=1}^6(5^i-1)).$$
If $n=6$, then $31,601\in \pi(G)\backslash \pi(\mbox{Aut}(S))$
and thus $31, 601\in \pi(N)$. Therefore $N$ contains an element
of order $31.601$, which is a contradiction as $31.601\notin
\omega(G)$. For case $n\geq 7$, since by Zsigmondy's Theorem $5_n,
5_{2n}>13$, a similar argument with the primitive prime divisors
$5_n, 5_{2n}\in \pi(G)$ also leads to a contradiction. Similarly,
$S$ can not be isomorphic to $S_4(3)$.
\begin{ca}
$q=7^n$, $n\geq 2$.
\end{ca}
In this case the possibilities for $S$ are: $ L_2(2^3)$,
$L_3(2^2)$, $S_6(2)$, $O_8^+(2)$, $G_2(3)$, $G_2(19)$, $U_3(3)$,
$U_3(5)$, $U_3(19)$, $U_4(3)$, $U_6(2)$ or $Sz(2^3)$. First of
all, since $G$ has no element of order $49$, $S\ncong G_2(19)$ or
$U_3(19)$. Next, we note that $\pi(S)\subset \pi(13!)$ and by
Lemma \ref{prim} we see that for every $n\geq 5$ there exists a
primitive prime divisor $7_n\geq 13$. Therefore for $n\geq 5$, as
previous cases a similar argument with the primitive prime
divisors $7_n$ and $7_{2n}$, leads to a contradiction.

If $n=4$, then $\mu(G)=\{2^5.3.5^2,7,2.1201\}$. In this case $S$
can only be  $L_2(2^3)$, $L_3(2^2)$, $S_6(2)$, $O_8^+(2)$,
$U_3(3)$, $U_3(5)$, or $U_4(3)$ by checking their prime
 divisors sets. On the other hand
, since the simple groups
 $L_2(2^3)$, $S_6(2)$, $O_8^+(2)$ and $U_4(3)$ contain an element
 of order 9 and $9\notin \omega(G)$, $S$ can only be $L_3(2^2)$,
 $U_3(3)$ or $U_3(5)$. Moreover, since there
does not exist any element of order 1201 in $\mbox{Aut}(S)$, 1201
divides the order of $N$. Without loss of generality we may
assume that $N\neq 1$ is an elementary Abelian 1201-group.
Because $A_4=2^2:3<A_6<L_3(2^2)$, $7:3<L_2(7)<U_3(3)$ and
$A_4=2^2:3<A_7<U_3(5)$, in all cases $S$ contains a Frobenius
group of shape $2^2:3$ or $7:3$ , and so $G$ contains an element
of order $1201.3$ by Lemma \ref{Maz3}, which is a contradiction.

If $n=3$, then $\mu(G)=\{2.3^2.19,7,2^3.43\}$. In this case, $S$
can only be $L_2(2^3)$ by checking their element orders sets.
 As $43\not \in \pi(\mbox{Aut}(S))$
we have $43\in \pi(N)$. Now, we may assume that $N\neq 1$ is an
elementary Abelian 43-group. Since $2^3:7<L_2(2^3)$ we get
$43.7\in \omega(G)$ by Lemma \ref{Maz3}, which is a contradiction.

If $n=2$, then $\mu(G)=\{2^4.3,7,2.5^2\}$. In this case, by
checking element orders $S$ can only be $L_3(2^2)$, $U_3(3)$ or
$U_3(5)$. If $S\cong L_3(2^2)$ or $U_3(5)$, then 5 divides the
order of $N$ since $25\notin \omega(\mbox{Aut}(S))$. Without loss
of generality we may assume that $N\neq 1$ is an elementary
Abelian 5-group. Since $S$ contains a Frobenius subgroup of shape
$2^2:3$ (in fact we have $A_4=2^2:3\leq A_6\leq L_3(4)$ and
$A_4=2^2:3\leq A_7\leq U_3(5)$), we get $5.3\in \omega(G)$ by
Lemma \ref{Maz3}, a contradiction. If $S\cong U_3(3)$, then $5\in
\pi(N)$, because $5\notin \pi (\mbox{Aut}(S))$. Again, since
$7:3\leq L_2(7)\leq U_3(3)$ we get $5.3\in \omega(G)$ by Lemma
\ref{Maz3}, a contradiction.
\begin{ca}
$q=13^n$, $n\geq 2$.
\end{ca}
In this case $S$ can only be isomorphic to one of  the following
 simple groups: $L_2(3^3)$, $L_2(5^2)$, $L_3(3)$, $L_4(3)$, $O_7(3)$,
$S_4(5)$, $S_6(3)$, $O_8^+(3)$, $G_2(2^2)$, $G_2(3)$, $F_4(2)$,
$U_3(2^2)$, $U_3(23)$, $Sz(2^3)$, $^3D_4(2)$, $^2E_6(2)$ or
$^2F_4(2)'$. Since $13^2\notin \omega(G)$, and $U_3(23)$ contains
an element of order $13^2$, $S\ncong U_3(23)$. Moreover, we have
$\pi(S)\subseteq \pi(19!)$. Now since, by Lemma \ref{prim}, for
every $n\geq 5$, there exists a primitive prime divisor $13_n>19$.
we can consider the primitive prime divisors $13_n$ and
$13_{2n}$, and we get a contradiction as before cases.
Henceforth, we may assume that $n\leq 4$.

If $n=4$, then $\mu(G)=\{2^4.3.5.7.17, 13, 2.14281\}$. In this
case, by comparing element orders, we conclude that $S$ can only
be $L_2(3^3)$, $L_2(5^2)$, $L_3(3)$, $L_4(3)$, $O_7(3)$,
$S_4(5)$, $O_8^+(3)$, $G_2(2^2)$, $G_2(3)$, $F_4(2)$, $U_3(2^2)$,
$Sz(2^3)$, $^3D_4(2)$, or $^2F_4(2)'$. In all above cases, except
$S\cong F_4(2)$ , since $17, 14281\not \in \pi(\mbox{Aut}(S))$,
we have $17, 14281\in \pi(N)$, and so $17.14281\in \omega(N)$,
which is a contradiction. If $S\cong F_4(2)$, then 14281 divides
the order of $N$, and since $S$ contains a Frobenius group
$2^2:3$(note that $2^2:3=A_4<S_{10}<S_8(2)<F_4(2)$), $G$ must
contain an element of order $14281.3$, by Lemma \ref{Maz3}, which
is not possible. If $n=3$, then $\mu(G)=\{2^2.3^2.61, 13,
2.7.157\}$. In this case we have $61, 157\not \in
\pi(\mbox{Aut}(S))$ and so $61, 157\in \pi(N)$, hence we get
$61.157\in \omega(N)\subset \omega(G)$, which is impossible.

\begin{ca}
$q=17^n$, $n\geq 2$.
\end{ca}
In this case $S$ can only be isomorphic to one of the following
simple groups: $L_2(2^4)$, $S_4(4)$, $S_8(2)$, $F_4(2)$,
$O_8^-(2)$, $O_{10}^-(2)$ or $^2E_6(2)$. First of all, since
$5\in \pi(S)$, we deduce $n\geq 4$. Moreover, we have
$\pi(S)\subseteq \pi(19!)$. From Lemma \ref{prim}, for every
$n\geq 4$, there exists a primitive prime divisor
 $17_n>19$.
Now, for the primitive prime divisors $17_n$ and $17_{2n}$, a
similar argument as  before leads to a contradiction.
\begin{ca}
$q=19^n$, $n\geq 2$.
\end{ca}
In this case $S$ can only be isomorphic to one of the following
simple groups:  $L_3(7)$, $U_3(2^3)$, $R(3^3)$ or $^2E_6(2)$.
Evidently $\pi(S)\subseteq \pi(37!)$. Since $7\in \pi(S)$, $3|n$.
If $n>7$, then by Lemma \ref{prim} there exists a primitive prime
divisor $19_n>37$. Now we consider the primes $19_n$ and
$19_{2n}$, and we get a contradiction as previous cases. If
$n=6$, then we have
$$\mu(G)=\{2^3.3^3.5.7.127, 19, 2.13^2.181.769\}.$$
In this case we consider the primes $127, 769 \in \pi(G)$, and we
obtain a contradiction as before. If $n=3$, then
$\mu(G)=\{2.3^3.127, 19, 2^2.5.7^3\}$. In this case $S$ can be
only $L_3(7)$ or $U_3(2^3)$, and since $5, 127\not \in
\pi(\mbox{Aut}(S))$, we get a contradiction.
\begin{ca}
$q=37^n$, $n\geq 2$.
\end{ca}
In this case $S$ can only be isomorphic to one of the following
simple groups: $U_3(11)$, $R(3^3)$ or $^2F_4(2^3)$. Evidently,
$\pi(S)\subseteq \{$2, 3, 5, 7, 11, 13, 19, 37, 73, 109$\}$. If
$n\geq 7$, then by Lemma \ref{prim} there exists a primitive
prime divisors $37_n>109$, and hence we consider the primes
$37_{n}, 37_{2n}\in \pi(G)$, and we get a contradiction as
before. Therefore we may assume that $n\leq 6$. Since

\begin{tabular}{ll}
$\pi(G)=\{2, 3, 5, 7, 13, 19, 31, 37, 43, 67, 137, 144061\}$, & $n=6$,\\
$\pi(G)=\{2, 3, 11, 19, 37, 41, 4271, 1824841\}$, & $n=5$,\\
$\pi(G)=\{2, 3, 5, 19, 37, 89, 137, 10529\}$, & $n=4$,\\
$\pi(G)=\{2, 3, 7, 19, 31, 37, 43, 67\}$, & $n=3$,\\
$\pi(G)=\{2, 3, 5, 19, 37, 137\}$, & $n=2$,\\
\end{tabular}

\noindent it is easy to see that $109\notin \pi(G)$, and so
$S\ncong {^2F_4(2^3)}$. Moreover, since $55\in\omega(U_3(11))
\backslash \omega(G))$, $S\ncong U_3(11)$. Finally, if $S\cong
R(3^3)$, since $13\in \pi(R(3^3))$, we must have $n=6$. Yet, in
this case, we can choose the primes $67, 144061\in
\pi(G)\backslash \pi(\mbox{Aut}(S))$, and we get a contradiction
as before (note that $67.144061\notin \omega(G)$).
\begin{ca}
$q=73^n$, $n\geq 2$.
\end{ca}
In this case $S$ can only be  $L_2(73^n)$, $L_3(2^3)$,
$S_6(2^3)$, $G_2(2^3)$, $G_2(3^2)$, $F_4(3)$, $E_6(2)$, $E_7(2)$,
$U_3(3^2)$ or $^3D_4(3).$ We assume that $S\ncong L_2(73^n)$. It
is not difficult to see that $\pi(S)\subseteq \pi(19!)\cup \{31,
41, 43, 73, 127\}$. Let $n\geq 5$. Then by Lemma \ref{prim}(6),
$73_n, 73_{2n}>127$. Evidently $73_n.73_{2n}\notin \omega(G)$, as
$73_{2n}\in \pi(73^n+1)$. On the other hand, since $73_n,
73_{2n}\notin \pi(\mbox{Aut}(S))$, $73_n, 73_{2n}\in \pi(N)$
which implies that $73_n, 73_{2n}\in \omega(N) \subseteq
\omega(G)$, a contradiction. Hence $n\leq 4$. Because

\begin{tabular}{ll}
$\omega(G)=\{2^5.3^2.5.13.37.41, 73, 2.14199121\}$, & $n=4$,\\
$\omega(G)=\{2^3.3^3.1801, 73, 2.7.37.751\}$, & $n=3$,\\
$\omega(G)=\{2^4.3^2.37, 73, 2.5.13.41\}$, & $n=2$,
\end{tabular}\\[0.2cm]
by checking the sets of element orders for each simple group, the
only possibility for $S$ is $U_3(3^2)$, when $n=4$. In this case,
we consider the primes 41 and 14199121 in $\pi(G)$. Since $41\in
\pi(73^4-1)$ and $14199121\in \pi(73^4+1)$, $41\nsim 14199121$
and also $41, 14199121 \notin \pi(\mbox{Aut}(S))$, which implies
that $41, 14199121\in \pi(N)$. Now by the nilpotency of $N$, we
obtain that $41.14199121\in \omega(N)\subset \omega(G)$, which is
a contradiction.

\begin{ca}
$q=109^n$, $n\geq 2$.
\end{ca}
The proof of this case follows immediately from Lemmas
\ref{item}(1) and \ref{pp}.

\begin{ca}
$q=(2^m+1)^n$, where $2^m+1$ is a prime and $n\geq 2$.
\end{ca}
In this case $S$ can only be isomorphic to: $L_2(2^m)$,
$S_a(2^b)$, $a=2^{c+1}$, $c\geq 1$, and $b=2^d$, $c+d=s$,
$F_4(2^e)$, $e\geq 1$, $4e=2^s$, $O_{2(m+1)}^-(2)$, $s>1$, or
$O_a^-(2^b)$, $a=2^{c+1}$, $c\geq 2$, and $b=2^d$, $c+d=s$.

If $S\cong L_2(2^m)$, then $\mu(\mbox{Aut}(S))=\{m, 2^m-1,
2^m+1\}=\{m, p-2, p\}$. First, assume that $n$ is odd. In this
case we have $(p-2,p^n-1)=1$, in fact if $(p-2,p^n-1)=d$ then $d$
divides $2^n-1$, and so $d\mid (p-2,2^n-1)=(2^m-1,
2^n-1)=2^{(m,n)}-1=1$. Now since $\pi(S)\subseteq \pi(G)$, it
follows that $\pi(p-2)\subset \pi(p^n+1)$. Moreover, it is easy
to see that $2^{m-1}+1$ divides $p^n+1$ and $(p-2,2^{m-1}+1)=3$.
Now we consider the primitive prime divisors
\begin{center}
$r:=p_n\in \pi(p^n-1)$  \ \ \  and  \ \ \  $s:=2_{2(m-1)}\in
\pi(2^{m-1}+1).$
\end{center}
Evidently $r, s\notin \pi(\mbox{Aut}(S))$, and so $r, s\in
\pi(N)$. From the nilpotency of $N$ it follows that $r\sim s$,
which is a contradiction. Next, we suppose that $n$ is even. In
this case we have $2^{m-1}+1$ divides $p^n-1$ and $(2^{m-1}+1,
p-2)=1$. Now, if $\pi(p-2)\subset \pi(p^n-1)$ then $(p-2,
p^n+1)=1$ and again we consider the following primitive prime
divisors
\begin{center}
$r:=p_{2n}\in \pi(p^n+1)$ \ \ \ and  \ \ \ $s:=2_{m-1}\in
\pi(2^{m-1}-1),$
\end{center}
and we get $r\sim s$, as before. But this a contradiction.
Therefore we must have $\pi(p-2)\subseteq \pi(p^n+1)$. Let $r\in
\pi(2^{m-1}+1)\subseteq \pi(p^n-1)$. Clearly $r\notin
\pi(\mbox{Aut}(S))$, hence $r\in \pi(N)$. Now since
$2^m:2^m-1\leq L_2(2^m)$, by Lemma \ref{Maz3} we deduce that
$r(2^m-1)\in \omega(G)$, which is a contradiction.

If $S\cong F_4(2^e)$, then the maximal odd factors set
$\psi(\mbox{Aut}(S))$ of $\mu(\mbox{Aut}(S))$ is equal to the
same set of $\mu(S)$ since $|\mbox{Out}(S)|=2^{e+1}$. From Lemma 9
we have
\begin{center}

$\psi(\mbox{Aut}(S))=\{q'^4-1, q'^4+1, q'^4-q'^2+1,
(q'-1)(q'^3+1), (q'+1)(q'^3-1)\},$

\end{center}
where $q'=2^e, e\geq 1$.

In this case $q'^4+1=p, q'^4-1=p-2.$ Since $G$ is an extension of
a $\pi_1(G)$-group $N$ by a group $G_1$, where $S\leq G_1 \leq
\mbox{Aut}(S)$, and $\mu(G)=\{p^n-1, p, p^n+1\}$, we may get a
contradiction dividing the two cases. If $n\geq 4,$ then the odd
number $\frac{1}{2}(p^n+1)$ and the odd factor of $p^n-1$ are all
greater than any number in $\psi(\mbox{Aut}(S)).$ Hence we have
$r,s$ such that
\begin{center}
$r\in \pi(p^n+1)$ \ \ \ and  \ \ \ $s\in \pi(p^n-1),$
\end{center}
and $r, s\notin \pi(\mbox{Aut}(S))$, so $r, s\in \pi(N)$. From the
nilpotency of $N$ it follows that $r\sim s$, which is a
contradiction. If $n=2,$ then we may infer that $(p-2, p^2+1)=5$
and $(p-2, p^2-1)=3.$ It is impossible. Also we may get a similar
contradiction if $n=3$.

If $S\cong S_a(2^b)$, then the maximal odd factors set
$\psi(\mbox{Aut}(S))$ of $\mu(\mbox{Aut}(S))$ is equal to the
same set of $\mu(S)$ since $|\mbox{Out}(S)|=b=2^d$. From
$\cite{kondratev}, \S3(3)$ we have
\begin{center}
$\{q'^{\frac{1}{2}(a)}-1, q'^{\frac{1}{2}(a)}+1\} \subseteq
\psi(\mbox{Aut}(S)),$
\end{center}
where $q'=2^b, b\geq 1$. In this case $q'^{\frac{1}{2}(a)}+1=p,$
and $q'^{\frac{1}{2}(a)}-1=p-2$, since the other numbers are not
primes in $\psi(\mbox{Aut}(S))$. The rest of proof is similar to
the case of $S\cong F_4(2^e)$ by comparing the two sets of
$\psi(\mbox{Aut}(S))$ and $\mu(G)$.

If $S\cong O_{2(m+1)}^-(2), m=2^s, s>1$, then the maximal odd
factors set $\psi(Aut(S))$ of $\mu(\mbox{Aut}(S))$ is equal to
the same set of $\mu(S)$ since $|\mbox{Out}(S)|=2$. From
$\cite{kondratev}, \S3(5)$ we have
\begin{center}
$\{q'^{m+1}+1, q'^m+1, q'^m-1\} \subseteq \psi(\mbox{Aut}(S)),$
\end{center}
where $q'=2$. In this case $q'^m+1=p,$ and $q'^m-1=p-2$. The rest
of proof is similar to the above cases.

If $S\cong O_a^-(2^b)$, $a=2^{c+1}$, $c\geq 2$, and $b=2^d$,
$c+d=s$, the proof is similar.
\begin{ca}
$q=97^n$ or $q=p^n$, where $p=2^\alpha 3^\beta +1>109$ is a prime,
$\beta \neq 0$ and $n\geq 2$.
\end{ca}
In this case $S$ is a simple $C_{pp}$-group, and from Table 1 and
Lemma \ref{item}, we obtain that $S\cong L_2(q)$.
\begin{st}
$N$ is a $2$-group.
\end{st}
Let $P/N$ be a Sylow $p$-subgroup of $S$ and $X/N$ be the
normalizer in $S$ of $P/N$. Then $X/N$ is a Frobenius group of
order $q(q-1)/2$, with cyclic complement of order $(q-1)/2$. Now,
by lemma \ref{Maz3}, we deduce that $N$ is a 2-group.
\begin{st}
$h(G)\in \{1,\infty \}$.
\end{st}
First suppose that $N=1$. In this case, we have $S=L_2(q)$,
$q=p^n$, $S\leq G \leq \mbox{Aut}(S)$. Denote the factor group
$G/S$ by $M$. Obviously, $M\leq \mbox{Out}(S)$. Therefore, every
element of $M$ is a product of a field automorphism $f$, whose
order is a divisor of $n$, and diagonal automorphism $d$ of order
dividing 2. Let $f\neq 1$ and $r$ be a prime dividing the order
of $f$. Without loss of generality, we may assume that $o(f)=r$.
Evidently, $r$ divides $n$, and we put $\bar{q}=p^{n/r}$. Denote
by $\varphi$ an automorphism of the field $\mathbb{F}_q$ inducing
$f$. Since $\varphi$ fixes a subfield $\mathbb{F}_{\bar{q}}$ of
$\mathbb{F}_q$, $f$ centralizes a subgroup $\overline{S}$ of $S$
isomorphic to $L_2(\bar{q})$. But then $G$ can not be a
$C_{pp}$-group, which is a contradiction. Thus $f=1$. Hence, we
have $M\leq \langle d\rangle$ and so $|G/S|\leq 2$. Therefore
$G\cong S$ or $G\cong \mbox{PGL}(2,q)$. From $q+1\in
\omega(\mbox{PGL}(2,q)) \backslash \omega(S)$, we have $G\cong
\mbox{PGL}(2,q)$. Thus, in this case $h(G)=1$. Next, suppose that
$N\neq 1$. Now, by Lemma
\ref{Maz1}, we get $h(G)=\infty$. The proof of Theorem \ref{th1} is complete. $\Box$\\[0.3cm]
\textsc{Proof of Theorem \ref{th2}.} Proof follows immediately
from
Lemma \ref{item}(2) and Lemma \ref{Maz1}. $\Box$\\[0.3cm]
\textsc{Proof of Theorem \ref{th3}.} Let $H$ be an extension of a
group of order 2 by $S_4$ such that a Sylow 2-subgroup of $H$ is a
quaternion group. Then $\mu(H)=\{6,8\}$. By Lemma 8 in \cite{m1},
there exists an extension $G$ of an elementary Abelian 7-group by
$H$, which is a Frobenius group. It follows that $\mu(G)=\{6, 7,
8\}$, and then Theorem \ref{th1} follows from Lemma \ref{Maz1}.
$\Box$
\begin{center}
{\Large \bf Acknowledgement}
\end{center}
The first author would like to thank the authority of K. N. Toosi
University of Technology very deeply for providing the necessary
facilities and their partial support to make the complete
performance of this project possible. The second author Supported
by the National Natural Science Foundation of Chaina (Grant No.
10171074).

\noindent {\sc A. R. Moghaddamfar} \\[0.2cm]
{\it Department of Mathematics, Faculty of Science,}\\
{\it  K. N. Toosi University of Technology,}\\ {\it P. O. Box 16315-1618, Tehran, Iran}\\
{\it E-mail:} {\tt moghadam@iust.ac.ir}\\[0.2cm]  {\sc W. J.
Shi}\\[0.2cm]
 {\it School of Mathematics, Soochow University,}
\\ {\it Suzhou 215006, People's Republic of China}\\
{\it E-mail:} {\tt wjshi@suda.edu.cn}
\end{document}